\documentclass[english,10pt]{article}
\usepackage[usenames,dvipsnames]{color}
\usepackage[colorlinks,linkcolor=OliveGreen,urlcolor=MidnightBlue,citecolor=Mahogany]{hyperref}%
\usepackage{hyperref}
\usepackage{mathrsfs}
\usepackage{bookmark}
 \usepackage[all]{hypcap}
\usepackage{amssymb, amsmath, latexsym}
\usepackage{layout}
\usepackage{geometry,graphicx,pict2e,authblk}
\usepackage[numbers,sort&compress]{natbib}
\usepackage{enumerate}
\usepackage{latexsym}
\usepackage{amssymb}
\usepackage{amsmath}
\usepackage{undertilde}
\usepackage{authblk}
\usepackage{multicol}
\usepackage{lipsum}
\usepackage[cp1252]{inputenc}
\usepackage[T1]{fontenc}

\newtheorem{Theo}{Theorem}[section]
\newtheorem{Lem}{Lemma}[section]

\newtheorem{Rem}{Remark}[section]

\newtheorem{Def}{Definition}[section]

\newcommand{\N}{{\rm  I~\hspace{-1.10ex}N}}

\newcommand{\ee}{\varepsilon}

\newcommand{\beq}{\begin{equation}}
\newcommand{\eeq}{\end{equation}}

\newcommand{\D}{\cal{D}}

\newcommand{\Ht}{\mathcal{H}_T}

\newcommand{\R}{\mathbb{R}}

\author{L. Boulanba  \thanks{Regional Center of the Trades of Education and Training, Souss-Massa. e-mail:
l.boulanba@gmail.com}\,\,\,\,\,\qquad  \qquad  \qquad \,M. Mellouk
\thanks{Universit\'{e} Paris Descartes,
e-mail:mohamed.mellouk@parisdescartes.fr}}

\title{Large deviations for a stochastic Cahn-Hilliard equation in H\"older norm}

\date{ \  \   }

\begin{document}

\maketitle 
\begin{abstract}
We consider a stochastic Cahn-Hilliard partial differential equation
driven by a space-time white noise.
 We prove the Large
Deviations Principle (LDP)  for the law of the solutions in the
H\"older norm. We use the weak convergence approach that reduces the
proof to establishing basic qualitative properties for controlled
analogues of the original stochastic system..
\end{abstract}
\section{Introduction}
  \hspace{0.3cm} In this paper we consider the following Stochastic Cahn-Hilliard
equation with multiplicative space-time white noise, indexed by
$\varepsilon>0$, given by
\begin{equation}
\left\{
\begin{array}{l}
\dfrac{\partial u^\varepsilon}{\partial t}(t,x)=- \Delta(\Delta
u^\varepsilon(t,x)-f( u^\varepsilon(t,x)))+ \sqrt{\varepsilon}\sigma
(u^\varepsilon(t,x)){\dot{W}}(t,x),\,\,\,(t,x)\in [0,T]\times
\mathcal{D},
\\
u^\ee(0,x)=u_0(x),\\
\dfrac{\partial u^\varepsilon}{\partial \nu}(t,x)=\dfrac{\partial
\Delta u^\varepsilon}{\partial \nu}(t,x)=0, \,\, \mbox{on}
\,\,\,[0,T]\times \partial \mathcal{D},
\end{array}%
\right.  \label{1}
\end{equation}%

\noindent where $T>0$,  $\mathcal{D}=[0,\pi]^d$ with $d=1,2,3$, $f$
is a polynomial of degree 3 with positive dominant coefficient such
as $f=F'$, where $F(u) = (1-u^2)^2$ is a double equal-well
potential. The noise diffusion coefficient $\sigma$ is a bounded and
Lipschitzian function, ${W}$ is a space-time Brownian sheet defined
on some filtered  probability space $(\Omega, {\cal F}, {\cal F}_t,
P)$, and $\nu$ is the outward normal vector. The initial condition
$u_0$ is a real-valued function satisfying some assumptions that
will be specified later.

 \noindent The deterministic Cahn-Hilliard equation (i.e., $\sigma \equiv 0$
in (\ref{1})) was introduced by Cahn and Hilliad in 1958 as a
mathematical model of spinodal decomposition for a binary allow in
order to determine the comprising species concentrations when the
separation phase take place, see \cite{Cahn-Hilliard}. In this
model, the function $f$ is the derivative of the homogeneous free
energy $F$ that was given in its original form by
$$
F(u) = -\frac{\theta_c}{2}u^{2}  +
\frac{\theta}{2}\left((1+u)\ln(1+u) - (1-u)\ln(1-u)\right), \qquad
-1 < u < 1,
 $$
 where $\theta_c$ and $\theta$ are respectively proportional to the critical
  and the quenching temperatures. One can see \cite{Debussche} where it is rigourously justified
  that $F$ can
 be replaced in some circumstances by a polynomial of even degree
 with a strictly positive dominant coefficient. For more details on physical
  aspects of this equation one can see e.g.
  (\cite{Elliott}, \cite{Langer}, \cite{Novick}).

\noindent Over the past three decades, different questions and
properties related to the equation (\ref{1}) have been the subject
of many works. Indeed, the existence of the solution and of its
density was established by Cardon-Weber in \cite{Cardon1}, a support
theorem was showed in \cite{BoShiWa10} by Bo and its co-authors, a
Freidin-Wentzell large deviations principle was obtained by Shi et
al. in \cite{S-T-W-09}. Recently, Antonopoulou et al.
\cite{Antonopoulou-Karali-Millet16} have attempted to go beyond
bounded coefficient noisy term in order to improve the results of
\cite{Cardon1}. Also, the Cahn-Hilliard equation driven by non
Gaussian perturbations was studied in a multitude of setting, we can
cite e.g. \cite{BoJiWa08} and \cite{JiShiWa14}.

\noindent Inspired by the pioneering works \cite{V} and \cite{F-W2}
 on large deviations for diffusion stochastic processes, growing
interest has been paid to this topic during the last three decades.
This was thanks to its various applications in many scientific
areas. Also, its nonlinear character and its connection with several
mathematic theories make it an active field of theoretical
researches. And besides a considerable literature about the large
deviations for stochastic differential equations (SDEs), this aspect
has been investigated for the most popular stochastic partial
differential equations (SPDEs) and we here cite e.g. \cite{MC-SA-03}
for the stochastic heat equation, \cite{Cardon2} for the stochastic
Burgers equation, \cite{Chenal-Millet97} for the stochastic wave
equation of degree two and \cite{Sowers92} for a reaction diffusion
equation with non- Gaussian perturbation. Note that in all these
works authors used the classical approach of Freidlin and Wentzell
that was developed essentially in \cite{Azencott80},
\cite{Priouret82} and \cite{Baldi-Cheleyat-Maurel88}. For a complete
and deep exposition of the topic of large deviations theory we refer
to \cite{D-Z}.

\noindent Recently, the weak convergence approach introduced by
Ellis and Dupuis in \cite{Ellis-Dupuis97} and developed in
\cite{Boue-Dupuis98}, \cite{B-D} and \cite{B-D-M} have gave a new
impetus to the study of large deviations both to investigate new
random dynamic systems or to revisit and improve anterior results of
the point of relaxing assumptions or simplifying the proof. And
taking advantage of this approach, many works on various SPDEs has
been appeared in last few years. See for a short list e.g.
(\cite{B-C-D}, \cite{Duan-Millet09}, \cite{El Mellali-Mellouk-16},
\cite{Seta14}, \cite{OrtSanz10}, \cite{Sritharan-Sundar06}). The
present paper fits into this optic.

\noindent It is worth mentioning that the weak convergence approach
consists to use a Laplace principle and some variational
representations for exponential functionals of infinite dimensional
Brownian motion. The proofs are based on showing qualitative
properties for controlled versions of the origin processes. This
fact unable one to avoid well known difficulties of the classical
approach when one wants establish exponential estimates that use
approximation and discretization procedure.

\noindent In this work we show a large deviations principle for the
stochastic Cahn-Hilliard equation in the H\"{o}der norm. Thereby, we
improve the result of \cite{S-T-W-09} that was given in terms of the
uniform convergence topology. Moreover, our proofs are technically
less demanding.

 \noindent The present paper is organized as
follows. Coming section contains basic backgrounds of large
deviations theory and well known results about the solution of the
equation (\ref{1}). Section 3 gives the general framework of our
work. In the last section, we announce and prove our main result.

\section{Preliminaries and main assumptions}

In this section we present some assumptions, preliminaries  and  standard
 definitions which are needed for  the formulation of the problem.\\
\subsection{Large deviations}
For a family of random variables  $\{X^\ee;\, \ee>0\}$ defined on a
probability space $(\Omega, \mathcal{F}, P)$ and taking values in a
Polish space ${\cal E}$, the LDP is concerned with events $A$ for
which probabilities $P(X^\ee \in A)$ converges to zero exponentially
fast as $\ee \rightarrow 0$. The exponential decay rate of such
probabilities are typically expressed in terms of a rate function
$I$ mapping ${\cal{E}}$ into $[0, \infty]$. 
 \begin{Def}
 The family of random variables $\{X^\ee;\, \ee>0\}$ is said to satisfy the
  LDP with the good rate function (or action functional)
 $I: {\cal{E}} \rightarrow  [0, \infty]$,  on ${\cal E}$, if
 \begin{enumerate}
 \item For each $M<\infty$ the level set $\{x\in  {\cal{E}};\,  I(x) \leq M\}$
 is a compact subset of  $\cal{E}$.
 \item \it{Large deviation upper bound:} for any closed subset $F $ of $\cal{E}$
$$\limsup_{\ee\rightarrow0^+}\ee\log P(X^{\ee}\in F)\leq -I(F).$$

\item \it{Large deviation lower bound:} for any open subset $O$ of  $\cal{E}$
$$\liminf_{\ee\rightarrow0^+}\ee\log P(X^{\ee}\in O)\geq -I(O).$$

\end{enumerate}
Where, for $A \subset\cal{E}$, we define $I(A)=\inf_{x \in A}I(x)$.
 \end{Def}
The Freidlin-Wentzell theory \cite{F-W} describes the path
asymptotics, as $\ee \rightarrow 0$, of probabilities of the large
deviations of the solutions of small noise finite dimensional SDEs,
away from its law if large number limite. For the case where the
driving
 brownian motion is infinite dimensional, that covers the SPDEs,
 Budhiraja et al. \cite{B-D-M} use certain variational representations to give a framework for proving large deviations
  for a variety of infinite dimensional systems. \\

\noindent In a many problems  one is interested in obtaining
exponential  estimates on functions which are more general than
indicator functions of closed or open sets. This leads to the study
of the, so called, Laplace principle.
\begin{Def}(Laplace principle)
 The family of random variables $\{X^\ee;\, \ee>0\}$ defined on the Polish space ${\cal E}$,
 is said to satisfy the Laplace principle with
rate function $I$
 if for any bounded continuous function $h:{\cal E} \rightarrow \R$,
{$$\lim_{\ee\rightarrow 0}\ee\log \mathbf{E}\left(
\exp\left[-\frac{1}{\ee}h(X^\ee)\right]\right)= -\inf_{f\in {\cal
E}}\{h(f)+I(f)\},$$}
where $\mathbf{E}$ is the expectation with respect to $P$.
\end{Def}

\noindent In \cite{V} and \cite{Br}, Varadhan and Bryc established
an equivalence between LDP and Laplace principle (LP) on a Polish
space.  In a view of this equivalence, the rest of this paper will
be concerned with the study of the Laplace principle.

\subsection{Assumptions and mild solution}
Letting $q\geq 1$, we define $\|\cdot \|_q$ as the usual norm in
$L^q(\D)$. Assume that: \vspace{0.3cm}
\begin{itemize}
 \item[(H1)] $f$ is a polynomial function of degree 3 with positive dominant
 coefficient. \vspace{0.3cm}
 \item[(H2)] $\sigma: $ is a bounded  and Lipschitz function.  \vspace{0.3cm}
 \item[(H3)] $u_0 \in L^p(\D)$ (for some $p\geq 4$) is continuous on
 $\D$.  \vspace{0.3cm}
 \item[(H3')] $u_0$ is an $\gamma$-H\"older continuous function on $\D$, $\gamma\in
]0,1]$.
\end{itemize}

\noindent Following the J. B. Walsh approach \cite{Wa}, a rigorous
meaning
 for solution of the equation (\ref{1}) can be given by means of the following definition.

\begin{Def}(Mild solution)
A jointly measurable and adapted process
$\{u(t,x);(t,x)\in [0.T]\times \D\}$  is called a mild solution of (\ref{1}) with initial condition $u_0$ if it satisfies, for each $t\geq 0$
and a.s. for almost all $x\in \D$ the following evolution equation:
\begin{eqnarray}
u^\ee(t,x) &=&\int_{\D}G_t(x,y)u_0(y)dy+\sqrt{\ee}\int_{0}^{t}\int_{\D}G_{t-s}(x,y)
\sigma(u^\ee(s,y))W(ds,dy) \notag \\
&&+\,\int_{0}^{t}\int_{\D}\Delta G_{t-s}(x,y)f(u^\ee(s,y))dsdy,
\label{3}
\end{eqnarray}%
where $G_{t}(\cdot,\cdot)$ denotes the Green kernel corresponding to
the operator $\frac{\partial}{\partial t}+\Delta^2$ with the Neumann
boundary conditions. Note that some useful estimates concerning
$G_{t}(\cdot,\cdot)$ are given in \cite{Cardon1}.
\end{Def}
The following result of C. Cadon-Weber (\cite{Cardon1}, Theorem 1.3)
asserts the existence and uniqueness of a solution to (\ref{1}).
\begin{Theo}(Existence, uniqueness and the regularity of solutions)
Under the assumptions (H1)-(H3),  there exists a unique  solution (in the Walsh's sense) of the equation \eqref{1}  which satisfies
\begin{equation}
\label{normelp} \sup_{0\le t \le T}\left( \mathbf{E}\|u(t,
\cdot)\|^{q}_{p}\right)^{1/q}< \infty,
\end{equation}
for $q\geq p$ if $d\in\{1,2\}$ and for $p\leq q <
\frac{6p}{(6-p)^{+}}$ if $d=3$. Moreover, under (H1)- (H3'),
(\cite{Cardon1}, Theorem 4.1) gives the a.s. $\beta$--H\"{o}lder
continuous property for the trajectories of the solution with $\beta
\leq \frac{\gamma}{4}$ and $\beta < \frac{1}{2}(1 - \frac{d}{4})$.
\end{Theo}
\begin{Theo}{(The  solution mapping of equation (\ref{1}))}.
\label{12a}
Assuming (H1)-(H3').  Let  $\alpha \leq
\frac{\gamma}{4}\wedge \frac{1}{2}(1 - \frac{d}{4})$ and set ${\cal E}_0=L^p({\D}) \cap {\cal C}^\gamma({\D})$, for some $p\geq 4 , \gamma \in ]0,1]$. There
exists a measurable function
$${\cal G}^\ee: {\cal E}_0 \times {\cal C}([0,T]\times {\D}:\R)\rightarrow C^{\alpha}([0,T], L^{p}(\D)),$$
 such that, for any filtered probability space $(\Omega, {\cal F}, {\cal F}_t, P)$ with a Brownian sheet $W$ as above and $u_0 \in {\cal E}_0$,
$u^{\ee}= {\cal G}^\ee(u_0, \sqrt{\ee}W), $
is the unique mild solution of (\ref{1}) (with initial condition $u_0$) and satisfies (\ref{normelp}).
\end{Theo}


 \hspace{0.2cm} In order to be able to apply the weak convergence approach for large deviations theory,
 we need a Polish space setting
carrying the probability laws of the family
$\{u^{\ee}(t,x);\,\ee\in(0,1],\,(t,x)\in[0,T]\times \mathcal{D}\}$.
And regarding the H\"{o}lder property of $u$ we introduce the space
$C^{\alpha}([0,T], L^{p}(\mathcal{D}))$ endowed with the following
norm
\begin{equation}
\|u\|_{\alpha, p} = \sup_{t\in [0,T]}\|u(t,\cdot)\|_{p} + \sup_{t\neq
t'\\
t, t' \in [0,T]}\frac{\|u(t,\cdot) - u(t',\cdot)\|_{p}}{|t-t'|^{\alpha}},
\end{equation}
for $p\geq 4$ and $\alpha \in ]0, 1[$. And because our setting
requires a Polish space state, we recall that if $\alpha' < \alpha
$, then the separable space $C^{\alpha',0}([0,T],
L^{p}(\mathcal{D}))$of functions belonging to $C^{\alpha'}([0,T],
L^{p}(\mathcal{D}))$ and satisfying
\begin{equation*}
\lim_{\delta \longrightarrow 0}\sup_{0 < |t-t'| < \delta, t\neq
t'}\frac{\|u(t,\cdot) - u(t',\cdot)\|_{p}}{|t-t'|^{\alpha'}} =0
\end{equation*}
 is a polish space containing $C^{\alpha}([0,T], L^{p}(\mathcal{D}))$. Hence,
 we can restrict ourselves in all the sequel to the space $\mathcal{E}^{\alpha}: =
C^{\alpha,0}([0,T], L^{p}(\mathcal{D}))$ for $\alpha <
\frac{\gamma}{4}\wedge \frac{1}{2}(1 - \frac{d}{4})$.

\section{Framework for the Laplace Principle}

Here, we first review an important result given by Budhiraja  et al.
\cite{B-D-M}, and which ensures the obtention of Laplace principle.

\subsection{Laplace principle of functionals of Brownian sheet.}


Consider the filtered probability space $(\Omega, {\cal F}, {\cal
F}_t, P)$ defined in the introduction, and let  ${\cal E}_0$ and
${\cal E}$ be Polish spaces such that the initial condition $u_0$
takes values in a compact subspace of ${\cal E}_0$.

\noindent Moreover, let $\left\{{\cal G}^\ee: {\cal E}_0 \times
{\cal C} ([0,T]\times { \D};
\R) \rightarrow {\cal E}, \  \ee>0 \right\}$  be a family of measurable maps.\\

\noindent For $u_0 \in {\cal E}_0$, Define
\begin{equation}
\label{BB} X^{\ee, u_0}: = {\cal G}^{\ee}(u_0,\sqrt{\ee}W).
\end{equation}

\noindent In the sequel we will give sufficient conditions for the
Laplace principle for $X^{\ee, u_0}$ to hold uniformly in $u_0$ for compact subsets of ${\cal E}_0$.\\

\noindent For $N \in \N$, consider the following
$$S^N = \left\{\phi \in
L^{2}([0,T]\times { \D}) : \int_{0}^{T} \int_{\D}\phi^2(s,y) ds dy
\leq N\right\},$$
\noindent which is a compact metric space, equipped with the weak topology on $L^{2}([0,T]\times { \D})$.\\

\noindent Let ${\cal P}_2$ be the class of all  predictable processes $\phi$ such that $\int_{0}^{T} \int_{\D}\phi^2(s,y) ds dy < \infty, a.s.$ Also,  define
$${\cal P}_2^N= \left\{v(\omega) \in {\cal P}_2:  v(\omega) \in S^N, P- a.s\right\},$$
the space of controls.\\
\noindent The following condition is the standing assumption of Theorem \ref{48} which states the Laplace principle for the family $\{X^{\ee, u_0}\}_{\ee>0}$ defined by (\ref{BB}). For $u\in L^{2}([0,T]\times { \D})$, define
${\cal I}(u) \in  {\cal C} ([0,T]\times { \D}: \R) $ as
$${\cal I}(u)(t,x):=\int_0^t \int_0^x u(s,y) ds dy.$$

\noindent {\sc \bf Assumption} $(\mathcal{A})$: There exists a
measurable map ${\cal G}^0$: ${\cal E}_0 \times {\cal C}
([0,T]\times { \D}; \R) \rightarrow {\cal E}$ such that the
following hold:
\begin{itemize}
 \item[(A1)] For every $M<\infty$ and a compact set $K\subset  {\cal E}_0$, the set
 $$ \Gamma_{M,K}:=
\left\{ {\cal G}^0(u_0,{\cal I}(u)); \ u \in S^M, u_0 \in K\right\}$$ is
a compact subset of $\cal{E}$.
\item[(A2)] Consider $M<\infty$ and a family $\{v^{\ee} : \ee > 0\} \subset {{\cal P}_2}^M$, and $\{ u_0^\ee \} \subset {\cal{E}}_0$
 such that $v^{\ee}\rightarrow v $ and $u_0^{\ee}\rightarrow u_0$
in distribution (as $S^N$-valued random elements)  as $\ee \rightarrow 0$. Then
 $${\cal{G}}^{\ee}(u_0^\ee, \sqrt{\ee}W+{\cal I}(v^\ee)) \rightarrow  {\cal G}^0(u_0,{\cal I}(u)),$$
in distribution as $\ee \rightarrow 0$.
\end{itemize}

\vspace{0.2cm}
 \noindent For $h \in {\cal E}$, and $u_0 \in {\cal
E}_0$, define the rate function
\begin{equation}
\label{rate}
I_{u_0}(h):=\inf_{\{  v\in L^{2}([0,T]\times { \D}): h= {\cal G}^0(u_0, {\cal I}(v)) \}} \left\{  \frac{1}{2}   \int_0^T \int_{\D} v^2(s,y)dy ds   \right\},
\end{equation}
 where the infimum over an empty set is taken to be $\infty$.\\

\noindent The following theorem  is due to
 Budhiraja  et al.(
\cite{B-D-M}, Theorem 7) and states the Laplace principle for the
family $X^{\ee, u_0}$.

\begin{Theo} (Theorem 7 in \cite{B-D-M})\label{48}
Let ${\cal G}^0$: ${\cal E}_0 \times {\cal C} ([0,T]\times { \D};
\R) \rightarrow {\cal E}$ be a measurable map satisfying assumption
$(\mathcal{A})$. Suppose that for all $h \in {\cal E}$,  $ u_0
\rightarrow I_{u_0}(h)$ is a lower semi-continuous map from ${\cal
E}_0$ to $[0, \infty]$. Then for every  $u_0 \in {\cal E}_0$,
$I_{u_0}(h): {\cal E}\rightarrow [0, \infty]$, is a rate function on
${\cal E}$ and the family $\{I_{u_0}, u_0 \in {\cal E}\}$ of rate
functions has compact level sets on compacts. Furthermore, the
family $X^{\ee, u_0}$ satisfies the Laplace principle on ${\cal E}$
with rate function $I_{u_0}$ defined by (\ref{rate}), uniformly in
$u_0$ on compact subsets of ${\cal E}_0$.
\end {Theo}
\subsection{The controlled and limiting equations for the spde (\ref{1})}
In the context of the spde under our study, ${\cal E}_0= L^p({\D}) \cap {\cal C}^\gamma({\D})$,
for some $p\geq 4 , \gamma \in ]0,1]$ is the space of the initial condition, and
${\cal E}=\mathcal{E}^{\alpha}: = C^{\alpha,0}([0,T], L^{p}(\D))$ for $\alpha < \frac{\gamma}{4}\wedge
\frac{1}{2}(1 - \frac{d}{4})$ the space of solutions. \\

\noindent The solution map of equation (\ref{1}) is $u^\ee={\cal
G}^\ee(u_0, \sqrt{\ee}W)$. Then, for $v\in {\cal P}_2^N$,
$u^{\ee,v}:= {\cal G}^\ee(u_0, \sqrt{\ee}W+ {\cal I}(v)) $ is the
solution map of the stochastic controlled equation for the spde
(\ref{1}) :

$$
\dfrac{\partial u^{\varepsilon,v}}{\partial t}(t,x)=- \Delta(\Delta
u^{\varepsilon, v}(t,x)-f( u^{\varepsilon,v}(t,x)))+
\sqrt{\varepsilon}\sigma
(u^{\varepsilon,v}(t,x)){\displaystyle\frac{\partial ^2 W}{\partial
t \partial x}}+ \sigma(u^{\varepsilon,v}(t,x))v(t,x),
$$
\noindent whose mild solution is
\begin{eqnarray}
\label{7}
 u^{\ee,v}(t,x)&=&\int_{\D}G_t(x,y)u_0(y)dy+\sqrt{\ee}\int_{0}^{t}\int_{\D}G_{t-s}(x,y)
 \sigma (u^{\ee,v}(s,y))W(ds,dy)\nonumber\\
& & + \int_{0}^{t}\int_{\D}\Delta G_{t-s}(x,y)f(u^{\ee,v}(s,y))dsdy
 + \int_{0}^{t}\int_{\D}G_{t-s}(x,y)\sigma(u^{\ee,v}(s,y))v^{}(s,y)dsdy.
\end{eqnarray}

\noindent Also, define the map $ {\cal G}^0(u_0, {\cal I}(v)): =
u^{v}$, where $u^{v}$ is the solution of the following zero-noise
equation:
\begin{eqnarray}
u^{v}(t,x)&=&\int_{\D}G_t(x,y)u_0(y)dy+
\int_{0}^{t}\int_{\D}G_{t-s}(x,y)\sigma (u^{v}(s,y))
v(s,y) ds dy\notag \\
& & + \int_{0}^{t}\int_{\D}\Delta G_{t-s}(x,y)f(u^{v}(s,y))dsdy.
\label{6-1}
\end{eqnarray}%

\noindent The following theorem gives a statement of existence and
uniqueness for the solution of the stochastic controlled equation
given by \eqref{7}.

\begin{Theo}{(Existence and uniqueness of controlled process)}
\label{12a}
Assuming $(H1)-(H3)$.  Let ${\cal G}^\ee$ denote the solution mapping, and let
 $v\in {\cal P}_2^N$ for some $N \in \N$.  Define
$$u^{\ee,v}= {\cal G}^\ee(u_0, \sqrt{\ee}W+ {\cal I}(v)), $$
then $u^{\ee,v}$ is the unique solution of equation (\ref{7}), which satisfies
\begin{equation}\label{estm}
 \sup_{\ee\leq1}\sup_{v\in{\cal P}_2^N}\sup_{0\le t\le T}
\mathbf{E}\left(\left\|u^{\ee,v}(t,\cdot)\right\|_{p}^q\right)<\infty,
\end{equation}
for $q \geq p$ if $d=1,2$, and $p\leq q < \frac{6p}{(6-p)^{+}}$ in
the case  $d=3$.

\end{Theo}

\noindent\textit{Proof.}
 For $ v\in {\cal P}_2^N$, set
$$
dQ^{\ee, v}:=
\exp\left\{-\frac{1}{\sqrt{\ee}}\int_{0}^{t}\int_{D}v(s,y)W(ds,dy)-
\frac{1}{2\ee}\int_{0}^{t}\int_{D}{v}(s,y)^{2}dsdy\right\}dP.
$$
\noindent Since is defined by an exponential martingale, $Q^{\ee,
v}$ is a probability measure on $\Omega$. And by Girsanov theorem
the process
$$
\tilde{W}(dt,dx) = W(dt,dx) +
\frac{1}{\sqrt{\ee}}\int_{0}^{t}\int_{D}v(s,y)dsdy
$$ is a
space-time white noise on the space $\Omega$ under the probability
measure $Q^{\ee, v}.$ Rewriting (\ref{7}) using $\tilde{W}(dt,dx)$
we obtain (\ref{3}) with $\tilde{W}(dt,dx)$ in place of
${W}(dt,dx)$.  Let $\textbf{u}$ be the unique solution of (\ref{3})
with $\tilde{W}(dt,dx)$ on the space $(\Omega, \mathcal{F}, Q^{\ee,
v})$.  Then $\textbf{u}$ satisfies (\ref{7}), $Q^{\ee, v}$ a.s. And
by equivalence of probabilities, then $\textbf{u}$ satisfies
(\ref{7}), $P$
a.s.\\
\noindent  For the uniqueness, if $\textbf{u}_1$ and $\textbf{u}_2$
are  two solutions of (\ref{7}) on $(\Omega, \mathcal{F}, P)$, then
$\textbf{u}_1$ and $\textbf{u}_2$ are solutions of (\ref{3})
governed by $\tilde{W}(dt,dx)$ on $(\Omega, \mathcal{F}, Q^{\ee,
v})$. By the uniqueness of the solution of (\ref{3}), we obtain
$\textbf{u}_1 = \textbf{u}_2$, $Q^{\ee, v}$ a.s. And by equivalence
of probabilities, we obtain $\textbf{u}_1 = \textbf{u}_2$, $P$ a.s.

\noindent Concerning the estimate (\ref{estm}), it holds true for
the three first terms by the estimations (2.16), (2.17) and (2.35)
in \cite{Cardon1}. It remains to show it for the last term. Indeed,
by the estimation (1.11) in \cite{Cardon1}, there exists a constant
$c > 0$ such that
\begin{eqnarray*}
\left\|\int_{0}^{t}\int_{D}G_{t-s}(\cdot,y)\sigma
(u^{\ee,v}(s,y))v(s,y)dsdy\right\|_{p}  \leq c
\int_{0}^{t}(t-s)^{\frac{d}{4}(\frac{1}{r}-1)}\|\sigma
(u^{\ee,v}(s,\cdot))v(s,\cdot)\|_{\rho}ds,
\end{eqnarray*}
where $\rho\in[1, p]$ and $\frac{1}{r} = \frac{1}{p} -
\frac{1}{\rho} + 1$. Using the boundedness of $\sigma$, taking $\rho
=2$ in the last inequality and applying Cauchy Schwarz inequality we
get a.s.
\begin{eqnarray}\label{123}
\left\|\int_{0}^{t}\int_{D}G_{t-s}(\cdot,y)\sigma
(u^{\ee,v}(s,y))v(s,y)dsdy\right\|_{p} & \leq&
\frac{c}{\frac{d}{2}\left(\frac{1}{r}-1\right) + 1}
T^{\frac{d}{2}\left(\frac{1}{r}-1\right) + 1}\|v\|_{\Ht}\nonumber\\
& \leq & \frac{c}{\frac{d}{2}\left(\frac{1}{r}-1\right) +
1}T^{\frac{d}{2}\left(\frac{1}{r}-1\right) + 1}N.
\end{eqnarray}
Note that, with the condition $p\geq 4$, there exists $r$ satisfying
(\ref{123}) that can be taken in $[\frac{4}{3}, 3[$.  Then
\begin{eqnarray}
\mathbf{E}\left(\left\|\int_{0}^{t}\int_{D}G_{t-s}(\cdot,y)\sigma
(u^{\ee,v}(s,y))v(s,y)dsdy\right\|_{p}^{q}\right)  < \infty.
\end{eqnarray}
Hence, (\ref{estm}) holds.  \hfill $\square$

\begin{Rem}{(H\"{o}lder regularity of controlled and limiting processes)}
\label{12} Assuming $(H1)-(H3')$. Both processes $\{u^{\ee,
v}(t,\cdot) ;\,t\in[0,T]\}$ and $\{u^{v}(t,\cdot) ;\,t\in[0,T]\}$,
defined by (\ref{7}) and (\ref{6-1}) respectively, live in the space
$\mathcal{E}^{\alpha}$.
\end{Rem}

\noindent\textit{Proof.}

\noindent The H\"{o}lder regularity for these two processes can be
obtained by arguing as in the point $ii)$ of the proof of Theorem
\ref{main}. \hfill $\square$

\noindent For $h \in {\cal E}^\alpha$, and $u_0 \in {\cal E}_0$, define the rate function

\begin{equation}
\label{rate2}
I_{u_0}(h):=\inf_{v} \left\{  \frac{1}{2}   \int_0^T \int_{\D} v^2(s,y)dy ds   \right\},
\end{equation}
 where  the infimum is taken over all $ v\in L^{2}([0,T]\times { \D})$ such that

\begin{eqnarray}
h(t,x)&=&\int_{\D}G_t(x,y)u_0(y)dy+
\int_{0}^{t}\int_{\D}G_{t-s}(x,y)\sigma (h(s,y))
v(s,y) ds dy\notag \\
&+&\int_{0}^{t}\int_{\D}\Delta G_{t-s}(x,y)f(h(s,y))dsdy.
\label{6}
\end{eqnarray}%

\noindent Note that under assumptions (H1)-(H3), for every
$v\in {\cal P}_2^N$, the equation \eqref{6} admits a unique solution which
belongs to $C([0,T], L^{p}(D))$, and moreover
\begin{equation}\label{30}
\sup_{t\in[0,T]} \|u^{v}(t, \cdot)\|^{q}_{p} < \infty,
\end{equation}
for $q\geq p$ if $d\in\{1,2\}$ and for $p\leq q \leq
\frac{6p}{(6-p)^{+}}$ if $d=3$. The proof is omitted since is
similar to that in Theorem 3.1 of \cite{Cardon1} but by replacing
the stochastic integral by the integral containing $v$.

\section{The main result}
\noindent The main result of this paper is the following:
\begin{Theo}\label{main}
 Under the assumptions (H1)-(H3'), the  law of the solution $\{u^{\ee};\,\ee\in(0,1]\}$,
 defined by \eqref{3}, satisfies, on $\mathcal{E}^{\alpha}$, a large deviation principle
with the rate function $I_{u_0}$, defined by (\ref{rate2}),
uniformly for $u_0$ in compact subsets of ${\cal E}_0$.
\end{Theo}

\noindent In view of Theorem \ref{48}, to prove Theorem \ref{main} it
suffices to verify conditions (A1) and (A2).

\begin{Rem}
This result improve that of Shi and al. \cite{S-T-W-09} where the
LDP was established by the classical approach in the space $C([0,T];
L^{p}(\mathcal{D}))$ equipped with the topology of uniform
convergence.
\end{Rem}

\noindent\textit{Proof of Theorem \ref{main}}
 \hspace{0.3cm} As mentioned above, here we will show that the
conditions $(A1)$ and $(A2)$
 hold. In a first time we deal with $(A2)$.
  That is, we need to show that
for all $q\geq p$ we have
\begin{equation}
\|u^{\ee,v^{\ee}}(t) - u^{v}(t)\|_{\alpha, p}^{q} \longrightarrow 0
\ \text{in probability as } \ \ {\ee\longrightarrow 0}.
\end{equation}
\noindent To do it, we will use a localization argument introduced
in \cite{Cardon-Millet}. For $0\leq t \leq T$, $\varepsilon \in
]0,1]$ and  $M
>0$, define the following event
\begin{equation*}
A_{\ee}^{M}(t) = \{ w\in \Omega; \sup_{s\in[0,t]}\|u^{\ee,
v^{\ee}}(s)\|_{p} \vee \sup_{s\in[0,t]}\|u^{v}(s)\|_{p} \leq M \}.
\end{equation*}
Notice that $A_{\ee}^{M}(t)\in \mathcal{F}_t$, and set
\begin{equation*}
Y_{\ee}(t): = u^{\ee,v^{\ee}}(t) - u^{v}(t).
\end{equation*}
\noindent Owing to (\ref{estm}) and (\ref{30}), we have
$P(A_{\ee}^{M}(T)^{c}) \longrightarrow 0$ \ \text{as} $\ \
{\ee\longrightarrow 0}$ and $ M \longrightarrow \infty$. Then, by
using Lemma A.1 in \cite{Cardon-Millet}, it suffices to show

\vspace{0.4cm}

\begin{itemize}
\item[$i)$] for all $ t\in [0,T];$
\begin{equation}
 \lim_{\ee \longrightarrow
0}\mathbf{E}\left[\mathbf{1}_{A_{\ee}^{M}(t)}\|Y_{\ee}(t)\|_{p}^{q}\right]=0
\label{22}
\end{equation}

\item[$ii)$] there exists $\beta > 0$ such that for all $ t, t'\in
[0,T]$,
\begin{equation}
 \sup_{\ee \in
]0,1[}\mathbf{E}\left[\mathbf{1}_{A_{\ee}^{M}(T)}\|Y_{\ee}(t) -
Y_{\ee}(t')\|_{p}^{q}\right]\leq c |t-t'|^{\beta q}. \label{23}
\end{equation}
\end{itemize}

\noindent  To prove $\ref{22}$, we write
\begin{eqnarray}
 Y_{\ee}(t) \nonumber
 & = & \sqrt{\ee} \int_{0}^{t}\int_{D}G_{t-s}(\cdot,y)
\sigma(u^{\ee,v^{\ee}}(s,y))W(ds,dy) + \int_{0}^{t}\int_{D}\Delta
G_{t-s}(\cdot,y)\left[f(u^{\ee,v^{\ee}}(s,y)) -
f(u^{v}(s,y))\right]dsdy,\nonumber\\
& & +  \int_{0}^{t}\int_{D}G_{t-s}(\cdot,y) \sigma (u^{\ee,v^{\ee}}(s,y))
\left[v^\ee(s,y) - v(s,y)\right] dsdy\nonumber\\
& & + \int_{0}^{t}\int_{D}G_{t-s}(\cdot,y)\left[\sigma (u^{\ee,v^{\ee}}(s,y)) -
\sigma (u^{v}(s,y))\right]v(s,y) dsdy\nonumber\\
 &  & \equiv \sum_{i=1}^{4} J_{i}^{\ee}(t).\nonumber
\end{eqnarray}
\noindent Then
\begin{eqnarray}\label{28}
\mathbf{E}\left[\mathbf{1}_{A_{\ee}^{M}(t)}\|Y_{\ee}(t)\|_{p}^{q}\right]
& \leq & c \sum_{i=1}^{4}
\mathbf{E}\left(\mathbf{1}_{A_{\ee}^{M}(t)}\|J_{i}(t)\|_{p}^{q}\right)\nonumber\\
& \leq & c\sum_{i=1,  i\neq 2}^{4}
\mathbf{E}\left(\|J_{i}(t)\|_{p}^{q}\right) +
c\mathbf{E}\left(\mathbf{1}_{A_{\ee}^{M}(t)}\|J_{2}(t)\|_{p}^{q}\right).
\end{eqnarray}

\noindent For $J_{1}^{\ee}(t)$, first we apply H\"{o}lder inequality
and we get
\begin{equation}
\mathbf{E} \left(\|J_{1}^{\ee}(t)\|_{p}^{q}\right)  = \mathbf{E}
\left(\int_{D} |J_{1}^{\ee}(t,x)|^{p}dx\right)^{\frac{q}{p}} \leq c
\int_{D} \mathbf{E} |J_{1}^{\ee}(t,x)|^{q}dx.
\end{equation}
\noindent Later we use Burkholder inequality, the boundedness of
$\sigma$ and the estimation (\ref{square green})
\begin{eqnarray}
\mathbf{E}|J_{1}^{\ee}(t,x)|^{q} & \leq & c \ee^{\frac{q}{2}}
\left(\int_{0}^{t}
\int_{D}G_{t-s}^{2}(x,y)\sigma^{2}(u^{\ee,v^{\ee}}(s,y))dsdy\right)^{q/2}\nonumber\\
& \leq & c \ee^{\frac{q}{2}} \left(\int_{0}^{t}\int_{D}G_{t-s}^{2}(x,y)dsdy\right)^{q/2}\nonumber\\
& < & c \ee^{\frac{q}{2}}.
\end{eqnarray}
 Concerning
$J_{2}^{\ee}(t)$, using (3.16) in \cite{BoJiWa08} and H\"{o}lder
inequality we get for $1 \leq \rho \leq p$ and $1 < \gamma \leq q$
\begin{eqnarray}
\mathbf{E}\left(\mathbf{1}_{A_{\ee}^{M}(t)}\|J_{2}^{\ee}(t)\|_{p}^{q}\right)
& \leq & c
\mathbf{E}\left(\mathbf{1}_{A_{\ee}^{M}(t)}\int_{0}^{t}\|f(u^{\ee,v^{\ee}}(s,\cdot))
-
f(u^{v}(s,\cdot))\|_{\rho}^{\gamma}ds\right)^{\frac{q}{\gamma}},\nonumber\\
& \leq & c
\mathbf{E}\left(\int_{0}^{t}\mathbf{1}_{A_{\ee}^{M}(s)}\|f(u^{\ee,v^{\ee}}(s,\cdot))
- f(u^{v}(s,\cdot))\|_{\rho}^{q}ds\right).\nonumber
\end{eqnarray}
\noindent Note that, for the last inequality, taking into account
the facts that  $A_{\ee}^{M}(t)\in \mathcal{F}_t$ and
$A_{\ee}^{M}(t) \subset A_{\ee}^{M}(s) $ for $0\leq s \leq t$, we
have used the following upper estimate
\begin{equation*}
\left|\mathbf{1}_{A_{\ee}^{M}(t)}\int_{0}^{t}\int_{D}\phi(s,y)dyds\right|
\leq
\left|\int_{0}^{t}\int_{D}\mathbf{1}_{A_{\ee}^{M}(s)}\phi(s,y)dyds\right|,
\end{equation*}
 for a measurable function $\phi: \Omega\times
[0,T]\times D \longrightarrow \mathbb{R}$. One can see Remark 3.2.
in \cite{Cardon-Millet}. \vspace{0.3cm}

\noindent And, since $f$ is a polynomial function of degree 3, we
can write
\begin{eqnarray*}
\|f(u^{\ee,v^{\ee}}(s,\cdot)) - f(u^{v}(s,\cdot))\|_{\rho}^{q} &
\leq & c \left[ \|u^{\ee,v^{\ee}}(s,\cdot) -
u^{v}(s,\cdot)\|_{\rho}^{q} + \|u^{\ee,v^{\ee}}(s,\cdot)^{2} -
u^{v}(s,\cdot)^{2}\|_{\rho}^{q}\right.\\
 & & \left. \quad + \|(u^{\ee,v^{\ee}}(s,\cdot)^{3} - u^{v}(s,\cdot)^{3}\|_{\rho}^{q}
\right].
\end{eqnarray*}
Taking $\rho = \frac{p}{3}$, we have
\begin{eqnarray*}
\hspace{-6cm} \|u^{\ee,v^{\ee}}(s,\cdot) -
u^{v}(s,\cdot)\|_{\rho}^{q} \leq c \|u^{\ee,v^{\ee}}(s,\cdot) -
u^{v}(s,\cdot)\|_{p}^{q};
\end{eqnarray*}
\begin{eqnarray*}
\hspace{-0.6cm} \|u^{\ee,v^{\ee}}(s,\cdot)^{2} -
u^{v}(s,\cdot)^{2}\|_{\rho}^{q} \leq c \|u^{\ee,v^{\ee}}(s,\cdot) -
u^{v}(s,\cdot)\|_{p}^{q}\left(\|u^{\ee,v^{\ee}}(s,\cdot)\|_{p}^{q} +
\|u^{v}(s,\cdot)\|_{p}^{q}\right);
\end{eqnarray*}
\begin{eqnarray*}
\|u^{\ee,v^{\ee}}(s,\cdot)^{3} - u^{v}(s,\cdot)^{3}\|_{\rho}^{q}
& \leq & c \|u^{\ee,v^{\ee}}(s,\cdot) - u^{v}(s,\cdot)\|_{p}^{q} \\
& &  \times \left(\|u^{\ee,v^{\ee}}(s,\cdot)\|_{p}^{q} +
\|u^{\ee,v^{\ee}}(s,\cdot)\|_{p}^{q}\|u^{v}(s,\cdot)\|_{p}^{q} +
\|u^{v}(s,\cdot)\|_{p}^{q}\right).
\end{eqnarray*}
\noindent Then
\begin{eqnarray}
\mathbf{E}\left(\mathbf{1}_{A_{\ee}^{M}(t)}\|J_{2}^{\ee}(t)\|_{p}^{q}\right)
\leq c
\int_{0}^{T}\mathbf{E}\left(\mathbf{1}_{A_{\ee}^{M}(s)}\|Y_{\ee}(s)\|_{p}^{q}\right)ds.
\end{eqnarray}

\noindent For $J_{3}^{\ee}(t)$, H\"{o}lder inequality, the
boundedness of $\sigma$ and the Cauchy schwarz inequality yield
\begin{eqnarray}
\mathbf{E}\left(\|J_{3}^{\ee}(t)\|_{p}^{q}\right) & \leq & c
\int_{D} \mathbf{E}
|J_{3}(t,x)|^{q}dx\nonumber\\
 & \leq & c \int_{D}\mathbf{E}\left( \int_{0}^{t}\int_{D}G_{t-s}(x,y)\left|v^\ee(s,y)
  - v(s,y)\right| dsdy\right)^{q}dx\nonumber\\
& \leq & c \sup_{x\in D}\left(\int_{0}^{t}\int_{D} G_{t-s}^{2}(x,y)dy\right)^{q}
\mathbf{E}\left(\|v^\ee - v \|_{\Ht}^{q}\right)\nonumber\\
& \leq & c \mathbf{E}\left(\|v^\ee - v \|_{\Ht}^{q}\right).
\end{eqnarray}

\noindent For $J_{4}^{\ee}(t)$, the same arguments as before yield
that a.s we have
\begin{eqnarray} \label{29}
\|J_{4}^{\ee}(t)\|_{p}^{q} & \leq & \|v\|_{\Ht}^{q}
\left(\int_{D}\left( \int_{0}^{t}\int_{D}G_{t-s}^{2}(x,y)
[\sigma (u^{\ee,v^{\ee}}(s,y)) -   \sigma (u^{v}(s,y))]^{2} dsdy\right)^{\frac{p}{2}}dx\right)^{q/p}\nonumber\\
& \leq  & N^{q}  \left( \int_{D} \left[
\left(\int_{0}^{t}\int_{D}G_{s}^{2}(x,y)dsdy\right)^{\frac{p}{2}-1}
 \int_{0}^{t}\int_{D}G_{t-s}^{2}(x,y) |Y_{\ee}(s,y)|^{p} dsdy\right]dx \right)^{q/p} \nonumber\\
 & \leq  & c  \left(\int_{D}\left[\int_{0}^{t}\int_{D}G_{t-s}^{2}(x,y)|u^{\ee,v^{\ee}}(s,y)
  -   u^{v}(s,y)|^{p}
 dsdy\right]dx \right)^{q/p} \nonumber\\
 & \leq  & c  \left(\int_{0}^{t}\left(\int_{D}G_{t-s}^{2}(x,y)dx
 \right)\int_{D}|Y_{\ee}(s,y)|^{p} dy ds\right)^{q/p}\nonumber\\
 & \leq  & c  \left(\int_{0}^{t} (t-s)^{-\frac{d}{4}}\|Y_{\ee}(s,y)\|_{p}^{p} ds\right)^{q/p}\nonumber\\
  & \leq  & c  \left(\int_{0}^{t} (t-s)^{-\frac{d}{4}}\right)^{\frac{q}{p} -1}
  \int_{0}^{t} (t-s)^{-\frac{d}{4}}\|Y_{\ee}(s,y)\|_{p}^{q} ds\nonumber\\
 & \leq  & c
 \int_{0}^{t}(t-s)^{-\frac{d}{4}}\|Y_{\ee}(s)\|_{p}^{q}ds.
\end{eqnarray}
\noindent Hence, combining (\ref{28})-(\ref{29}) we obtain
\begin{eqnarray*}
\mathbf{E}\left(\mathbf{1}_{A_{\ee}^{M}(t)}\|Y_{\ee}(t)\|_{p}^{q}\right)
\leq c \left( \ee^{\frac{q}{2}} +  \mathbf{E} \left(\|v^\ee - v
\|_{\Ht}^{q}\right) + \int_{0}^{T}\left(1 +
(t-s)^{-\frac{d}{4}}\right)
\mathbf{E}\left(\mathbf{1}_{A_{\ee}^{M}(s)}\|Y_{\ee}(s)\|_{p}^{q}\right)ds
 \right).
\end{eqnarray*}
We obtain (\ref{22}) by applying a version of Gronwall lemma given
by (Lemma 15, \cite{Dalang99}).

\vspace{0.5cm}

\noindent To prove $(\ref{23})$, consider $t$; $t'\in [0,T]$ such
that $t < t'$. We
 have
 $$
Y_{\ee}(t) - Y_{\ee}(t') = (u^{^{\ee,v^{\ee}}}(t) -
u^{^{\ee,v^{\ee}}}(t')) - (u^{v}(t)  - u^{v}(t')).
 $$
 Then
$$
\mathbf{E}\left(\mathbf{1}_{A_{\epsilon}^{M}(T)}\|Y_{\ee}(t) -
Y_{\ee}(t')\|_{p}^{q}\right) \leq
\mathbf{E}\left(\mathbf{1}_{A_{\epsilon}^{M}(T)}\|u^{^{\ee,v^{\ee}}}(t)
- u^{^{\ee,v^{\ee}}}(t')\|_{p}^{q}\right) +
 P({A_{\epsilon}^{M}(T)})\|u^{v}(t)  - u^{v}(t')\|_{p}^{q}.
$$
At beginning we deal with the first term and we write (\ref{7}) as
$$
u^{\ee,v^{\ee}}(t, \cdot) :=  I_{0}^{\ee,v^{\ee}}(t,\cdot) +
I_{1}^{\ee,v^{\ee}}(t,\cdot) + I_{2}^{\ee,v^{\ee}}(t,\cdot) +
I_{3}^{\ee,v^{\ee}}(t,\cdot),
$$
where $I_{i}^{\ee,v^{\ee}}(t,x)$ stands for the $i$-th term in
(\ref{7}).  By Lemma 2.2. in \cite{Cardon1}, there exists a constant
$c>0$ such that
\begin{equation}\label{24}
\|I_{0}^{\ee,v^{\ee}}(t) - I_{0}^{\ee,v^{\ee}}(t')\|_{p}^{q} \leq c
|t -t'|^{\frac{q\gamma}{4}}.
\end{equation}
The same reference gives the existence of $\beta >0$ such that $0 <
\beta <\frac{1}{2}\left(1-\frac{d}{4}\right)$ such that
\begin{equation}\label{25}
\mathbf{E}\left(\|I_{1}^{\ee,v^{\ee}}(t) -
I_{1}^{\ee,v^{\ee}}(t')\|_{p}^{q}\right)\leq c\sqrt{\ee}|t -
t'|^{\beta q}.
\end{equation}
Concerning $I_{2}^{\ee,v^{\ee}}(t)$ we have
\begin{eqnarray*}
I_{2}^{\ee,v^{\ee}}(t') - I_{2}^{\ee,v^{\ee}}(t) & = &
\int_{0}^{t}\int_{D}\Delta
[G_{t'-s}(\cdot,y) - G_{t-s}(\cdot,y)]f(u^{\ee,v^{\ee}}(s,y))dsdy \\
& &  + \int_{t}^{t'}\int_{D}\Delta
G_{t'-s}(\cdot,y)f(u^{\ee,v^{\ee}}(s,y))dsdy\\
& \equiv &  I_{2,1}^{\ee,v^{\ee}}(t,t') +
I_{2,2}^{\ee,v^{\ee}}(t,t').
\end{eqnarray*}
\noindent By (1.12) in \cite{Cardon1}, the H\"{o}lder inequality and
the estimation (\ref{estm})  there exists $1\leq \rho \leq p$ and
$\kappa \in [0,1]$ such that
\begin{eqnarray}
\mathbf{E}\left(\|I_{2,2}^{\ee,v^{\ee}}(t,t')\|_{p}^{q}\right) &
\leq & c \mathbf{E}\left(\int_{0}^{t'-t}(t'-t-s)^{-\frac{1}{2} +
\frac{d}{4}(\kappa
-1)}\|f(u^{\ee,v^{\ee}}(t+s,\cdot))\|_{\rho}ds\right)^{q}\nonumber\\
& \leq & c |t' -t|^{q\left(\frac{1}{2} + \frac{d}{4}(\kappa
-1)\right)} \int_{0}^{t'-t}(t'-t-s)^{-\frac{1}{2} +
\frac{d}{4}(\kappa
-1)}\mathbf{E}\left(\|f(u^{\ee,v^{\ee}}(t+s,\cdot))\|_{\rho}^{q}\right)ds\nonumber\\
& \leq & c |t' -t|^{(q  +1)\left(\frac{1}{2} + \frac{d}{4}(\kappa
-1)\right)}.
\end{eqnarray}
\noindent Using (3.14) in \cite{BoJiWa08}, H\"{o}lder inequality and
the estimation (\ref{estm}), there exist $\theta \in \left]0,
\frac{1}{2} + \frac{d}{4}(\kappa -1)\right[$ and $\gamma\in
\left]\frac{1}{\frac{1}{2} + \frac{d}{4}(\kappa -1)-\theta}, q
\right[$ such that
\begin{eqnarray}
\mathbf{E}\left(\|I_{2,1}^{\ee,v^{\ee}}(t,t')\|_{p}^{q}\right) &
\leq & c
 |t' -t|^{\theta q}\mathbf{E}\left(\|f(u^{\ee,v^{\ee}}(\cdot,\ast))\|_{L^{\gamma}([0,T], L^{\rho}(D))}^{q}\right)\nonumber\\
& \leq & c |t' -t|^{\theta q}. \label{26}
\end{eqnarray}

\noindent Concerning $I_{3}^{\ee,v^{\ee}}(t)$ we have
\begin{eqnarray*}
I_{3}^{\ee,v^{\ee}}(t') - I_{3}^{\ee,v^{\ee}}(t) & = &
\int_{0}^{t}\int_{D}[G_{t'-s}(\cdot,y) - G_{t-s}(\cdot,y)]\sigma(u^{\ee,v^{\ee}}(s,y))v^{\ee}(s,y) dsdy \\
& &  + \int_{t}^{t'}\int_{D} G_{t'-s}(\cdot,y)\sigma(u^{\ee,v^{\ee}}(s,y))v^{\ee}(s,y)dsdy\\
& \equiv &  I_{3,1}^{\ee,v^{\ee}}(t,t') +
I_{3,2}^{\ee,v^{\ee}}(t,t').
\end{eqnarray*}
By Cauchy-Schwarz inequality, the fact that $\|v^{\ee}\|_{\Ht} \leq
N$ a.s. and by Lemma 1.8. in \cite{Cardon1} we obtain the existence
of $\eta > 0$ such that
\begin{eqnarray}\label{27}
\mathbf{E}\left(\|I_{3,i}^{\ee,v^{\ee}}(t,t')\|_{p}^{q}\right) &
\leq & c |t' -t|^{\eta q},
\end{eqnarray}
for $i = 1, 2$.

\noindent Therefore, by (\ref{24})--(\ref{27}) we obtain (\ref{23})
for the first term. And arguing similarly and using the estimation
(\ref{30}) we obtain (\ref{23}) for the second term. Hence, the
condition (A2) is checked.

\vspace{0.3cm} \noindent Concerning $(A1)$, it will be a consequence
of the continuity of the mapping $h: \mathcal{H}_{T}^{N}
\longrightarrow \mathcal{E}^{\alpha}$ with respect to the weak
topology. It consists to consider $v$, $(v_n)\subset
\mathcal{H}_{T}^{N}$ such that for any $g\in \mathcal{H}_{T}^{N}$,
$$
\lim_{n\longrightarrow +\infty} \langle v - v_n , g
\rangle_{\mathcal{H}_{T}^{N}} =0,
$$
and to prove
\begin{equation}
\lim_{n \longrightarrow +\infty} \|u^{v_n}  - u^{v}\|_{\alpha ,
p}=0.
\end{equation}
The proof will be omitted since we can proceed as for $(A2)$ and by
using the following estimate
\begin{equation}
\sup_{\|v\| \leq N}\sup_{t\in [0,T]}\|u^{v}(t)\|_{ p} < \infty,
\end{equation}
which follows from Lemma 3.1. in \cite{S-T-W-09}.

\vspace{0.3cm} \noindent Finally, since the conditions $(A1)$ and
$(A2)$ are held, the proof of Theorem \ref{main} is completed.\hfill
$\square$

\section{Appendix}
 \hspace{0.3cm} We recall here some useful results that we have used in the proofs of our result.
 The following lemma gives well-known estimates on space and time increments for the Green function $G$. For the proof, we refer
 to \cite{Cardon1}  .

 \begin{Lem}
 There exists positive constants $c$, $\gamma$ and $\gamma'$ satisfying
 $\gamma < 4-d$, $\gamma\leq 2$ and $\gamma'< 1-\frac{d}{4}$ such
 that for all $y,z \in \mathcal{D}$, $0\leq s < t \leq T$ and $0\leq h \leq t$ we have
 \begin{enumerate}
 \item
\begin{equation}
\int_{0}^{t}\int_{\mathcal{D}}|G_{r}(x,y) - G_{r}(x,z)|^{2}dxdr \leq
c |y-z|^{\gamma}
\end{equation}

 \item
\begin{equation}
\int_{0}^{t}\int_{\mathcal{D}}|G_{r + h}(x,y) - G_{r}(x,y)|^{2}dxdr
\leq c |h|^{\gamma'}
\end{equation}
 \item
\begin{equation}\label{square green}
\int_{s}^{t}\int_{\mathcal{D}}|G_{r}(x,y)|^{2}dxdr \leq c
|t-s|^{\gamma'}
\end{equation}
\end{enumerate}
\end{Lem}

\noindent The following lemma is a version of the
Garsia-Rademich-Rumsay lemma. For the proof, we refer to
\cite{Cardon-Millet} and references therein.

\begin{Lem}
Let  $\alpha \ ]0; 1]$ and  $1 < p \leq q$. Consider a sequence of
stochastic processes $\left(Y_n\right)_n$ which belong to
$C^{\alpha}([0,T]; L^{p}(\mathcal{D}))$, and a sequence of stopping
times $\left(\tau_n\right)$ such that
\begin{enumerate}
\item   for any $t\in [0,T]$,
$$
\displaystyle \lim_{n\longrightarrow
+\infty}\mathbf{E}\left(1_{\{t\leq
\tau_n\}}\|Y_{n}(t,\cdot)\|_{p}^{q}\right) = 0;
$$
\item there exists $\gamma > 0$ such that for any $(t,t')\in [0,T]$
$$
\displaystyle \sup_{n}\mathbf{E}\left(1_{\{t\vee t'\leq
\tau_n\}}\|Y_{n}(t,\cdot) - Y_{n}(t', \cdot)\|_{p}^{q}\right) \leq c
|t-t'|^{\gamma + d},
$$
\end{enumerate}
then, for any $1\leq r < q$ and any $\theta < \frac{\gamma}{q}$ one
has
$$
\displaystyle \lim_{n\longrightarrow
+\infty}\mathbf{E}\left(1_{\{t\leq
\tau_n\}}\|Y_{n}(t,\cdot)\|_{\theta, p, \tau_n}^{r}\right) = 0;
$$
where, for a stopping time $\tau$

\begin{equation*}
\|u\|_{\theta, p, \tau}: = \sup_{t\in [0,T\wedge\tau]}\|u(t)\|_{p} +
\sup_{t\neq
t'\\
t, t' \in [0,T\wedge\tau]}\frac{\|u(t) -
u(t')\|_{p}}{|t-t'|^{\theta}}.
\end{equation*}
\end{Lem}

\vspace*{1.5cm}
{\scshape
\begin{flushright}
\begin{tabular}{l}
D\'{e}partement de Math\'{e}matiques\\
Centre R\'{e}gional des M\'{e}tiers de l'\'{E}ducation et de la Formation\\
80 000 Agadir \\
Maroc\\
{\upshape e-mail: \href{l.boulanba@gmail.com}{\nolinkurl{l.boulanba@gmail.com}}},\\
\\
MAP5, CNRS UMR 8145 \\
Universit\'{e} Paris Descartes \\
45, rue des Saints-P\`eres \\
75270 Paris Cedex 6 \\
France \\
{\upshape e-mail: \href{mailto: mohamed.mellouk@parisdescartes.fr}{\nolinkurl{mohamed.mellouk@parisdescartes.fr}}}\\
\end{tabular}
\end{flushright}
}

\end{document}